\documentclass[preprint]{article}

\usepackage{amsmath}
\usepackage{amsthm}
\usepackage[all]{xy}

\newcommand{\NN}{\mathbf{N}}
\newcommand{\ZZ}{\mathbf{Z}}
\newcommand{\QQ}{\mathbf{Q}}
\newcommand{\CC}{\mathbf{C}}
\newcommand{\CALD}{\mathcal{D}}

\newtheorem{theorem}{Theorem}
\newtheorem{proposition}{Proposition}
\newtheorem{lemma}{Lemma}
\newtheorem{corollary}{Corollary}

\begin{document}

\title{Lambda-conductors for group rings}

\author{F. J.-B. J. Clauwens}

\maketitle

\section{Introduction.}

This paper is part of a  project which aims to
provide a method for computing the Nil groups of
the group rings of finite abelian groups,
by refining some of the techniques used in \cite{ads} and \cite{ados}
in such a way that the allowed coefficient rings include polynomial rings.
For the refinement of the $p$-adic logarithm discussed in 
\cite{flog} and \cite{fexp} it is assumed that the rings involved
have a structure of $\lambda$-ring;
we refer to these papers for generalities about $\lambda$-rings.
Thus it is useful to extend as much as possible of the  other techniques
to the context of $\lambda$-rings.
In this paper we investigate how to describe a group ring
of a finite abelian group as a pull back of a diagram of rings
which are more accessible to calculations in algebraic K-theory.

Let be given a commutative ring $S$ and subring $R$.
For each ideal $I$ of $S$ which is contained in $R$ one has a cartesian square
\begin{equation*}
\xymatrix{
R\ar[r]\ar[d]&S\ar[d]\\R/I\ar[r]&S/I
}
\end{equation*}
thus describing $R$ as a pull back of rings for which the $K$-theory is  hopefully
better understood.
By taking for $I$ the sum of all such ideals one finds a diagram
where the rings on the bottom row are as small as possible.

We modify this construction by assuming 
that $R$ has a structure of $\lambda$-ring and considering only ideals 
$I$ stable under the $\lambda$-operations.
We call the resulting ideal the $\lambda$-conductor of $S$ into $R$.
In particular we are interested in the case that $R$ is the group ring
$\ZZ[G]$ of a finite abelian group, and $S$ is its normal closure
in $R\otimes\QQ$, which splits as a direct sum of rings $S_i=\ZZ[\chi_i]$
associated to equivalence classes of characters $\chi_i\colon G\to\CC$.

In this situation $R$ is a $\lambda$-ring such that
$\psi^n(g)=g^n$ for every $n\in\NN$ and $g\in G$.
In general $S$ is not stable under the 
$\lambda$-operations on $R\otimes\QQ$,
but it is stable under the associated Adams operations $\psi^n$
since they are ring homomorphisms.

We will prove that in case $G$ is a primary group
its $\lambda$-conductor is precisely the intersection of 
the classical conductor and the augmentation ideal.
We do this by exhibiting generators of the classical
conductor and examining their behavior 
under the fundamental $\lambda$-operations.

\section{The primary case}

Throughout this section  $G$ is a group of order $n=p^e$, 
where $p$ is prime.
We consider representations $\rho\colon G\to \CC^*$.
We say that $\rho$ is of level $k$ if the image of $\rho$ has $p^k$ elements.

Two representations $\tau_1$ and $\tau_2$ are called equivalent
if they have the same kernel.
That means that there must be some $m\in\ZZ$ prime to $p$ such that
$\tau_2(x)=\tau_1(x)^m$ for all $x\in G$.
Obviously equivalent representations have the same level.

Given a representation $\tau$ of level $k>0$ one gets a representation
$\psi\tau$ of level $k-1$ by the formula $(\psi\tau)(x)=\tau(x^p)$ for $x\in G$.
If $\psi\tau_1$ and $\psi\tau_2$ are equivalent  then we may replace
$\tau_2$ by an equivalent representation $\tau'_2$
so that $\psi\tau'_2=\psi\tau_1$.
So we may choose a representation in each class in such 
a way that $\psi\tau$ and $\rho$ coincide if they are equivalent.

Let $\rho$ be a representation of level $k>0$ 
and write $\omega=\exp(2\pi i/p)$.
We define an element $b_\rho\in\ZZ[G]$ by the formula
\begin{equation*}
b_\rho=\sum_{\rho x=1}x-\sum_{\rho\xi=\omega}\xi.
\end{equation*}
If we choose $y_\rho\in G$ such that $\rho(y)=\omega$ then we get
\begin{equation*}
b_\rho=\left(\sum_{\rho x=1} x\right)\left(1-y_\rho\right)
\end{equation*}
The only representation of level $0$ is the trivial representation,
which we denote by $1$, and it gives rise to 
$b_1=\sum_{x\in G}x$.

\begin{proposition}
\label{bortho}
If  $\rho$ and $\tau$ are not equivalent then $b_\rho b_\tau=0$.
Furthermore $b_\rho^2=p^{e-k}\left(1-y_\rho\right)b_\rho$
for $\rho$ of level $k>0$.
\end{proposition}

\begin{proof}
If $\ker(\rho)\not=\ker(\tau)$ we may assume that there is
$g\in G$ with $\rho(g)=1$ but $\tau(g)=\omega$.
If $g$ has order $m$ then $\sum_{\rho(x)=1} x$ and thus $b_\rho$
is a multiple of $\sum_{j=0}^{m-1} g^j$,
whereas $b_\tau$ is a multiple of $1-g$ .
The product of these two factors is $0$.

The second part follows from
$(\sum_{\rho(\xi)=1}\xi)(\sum_{\rho(x)=1}x)=p^{n-k}\sum_{\rho(x)=1}x$
which is true because each $\xi$ gives the same contribution 
and there are $p^{e-k}$ of them.
\end{proof}

Every representation $\rho$ of level $k$ gives rise to a homomorphism $j_\rho$
from $\ZZ[G]$ to $S_\rho=\ZZ[\omega_k]$,
where $\omega_k=\exp(2\pi i/p^k)$.\begin{proposition}
\label{jb}
If  $\rho$ and $\tau$ are not equivalent then $j_\rho(b_\tau)=0$.
Furthermore  $j_1(b_1)=p^e$, and $j_\rho(b_\rho)=p^{e-k}\left(1-\omega\right)$
if $\rho$ is of level $k>0$,.
\end{proposition}

\begin{proof}
The second part  is obvious since every $x$ in the definition of $b_\rho$
maps to $1$, and $y_\rho$ maps to $\omega$.
The first part follows since $j_\rho(b_\rho)j_\rho(b_\tau)=0$
by Proposition \ref{bortho} and $S_\rho$ is a domain.
\end{proof}

It is well known that the maps $j_\rho$ (one from each equivalence class)
combine to an embedding from $R$ into its integral closure
$S=\oplus_\rho S_\rho$.
\begin{proposition}
\label{conduc}
The $b_\rho$ generate the conductor ideal $I$ of $S$ into $R$.
\end{proposition}

\begin{proof}
By the theorem of Jacobinski (Theorem 27.8 in \cite{curt})
the conductor is $\oplus_\rho n\CALD_\rho^{-1}\subset\oplus_\rho S_\rho=S$,
where $\CALD_\rho^{-1}$ is the lattice in $\QQ[\omega_k]$
dual to $\ZZ[\omega_k]$ under  the trace form.
It is a simple exercise that this fractional ideal is in fact
generated by $p^{-k}(1-\omega)$,
which means that $n\CALD_\rho^{-1}$ is 
just $j(b_\rho)S=j(B_\rho)S_\rho=j(b_\rho R)$.
\end{proof}

We remind the reader that in particular $nS$ is contained in the conductor.

\bigbreak

The $\lambda$-conductor $I_\lambda$ from $S$ into $R$ is
defined as the largest ideal of $S$ contained in $R$
which is stable under the fundamental $\lambda$-operations $\theta^\ell$.
It is of course a subset of the largest ideal of $S$ contained in $R$,
which is the ordinary conductor $I$ described above.
Thus we have to investigate the behaviour of the operations
$\theta^\ell$ on the generators $b_\rho$.

\begin{lemma}
If $\rho$ is of level $k>0$ and
there is no $\tau$ with $\psi\tau=\rho$
then $\psi^pb_\rho=0$.
\end{lemma}

\begin{proof}
Write $G$ as a direct product of cyclic groups,
with generators $g_i$.
If the order of $\rho(g_i)$ is strictly smaller than the order of $g_i$
for all $i$ then one can find 
a suitable $\tau$ by taking for each $\tau(g_i)$
a $p$-th root of $\rho(g_i)$.
If however the orders are the same for some $i$ then 
there is certainly some $h\in G$ such that
$h^p=1$ and $\rho(h)=\omega$.
By definition of $b_\rho$ we have
\begin{equation*}
\psi^p b_\rho=\sum_{\rho\xi=1}\xi^p-\sum_{\rho\eta=\omega}\eta^p
\end{equation*}
Here the term in the second sum associated to $\eta=h\xi$
cancels the term in the first sum associated to $\xi$.
\end{proof}

\begin{proposition}
If $\rho$ is of level $k>0$ then
\begin{equation*}
\psi^p(b_\rho)=\sum_{\psi\tau=\rho} p b_\tau
\end{equation*}
\end{proposition}

\begin{proof}
By definition we have
\begin{equation*}
\sum_{\psi\tau=\rho}b_\tau=
\sum_{\psi\tau=\rho}\;\sum_{\tau x=1}x
-\sum_{\psi\tau'=\rho}\;\sum_{\tau' x=\omega}x
\end{equation*}
We claim that all terms with $x\not\in G^p$ cancel.
To prove this assume that the class of $x$ in
$G/G^p$ is nontrivial.
Then there exists a homomorphism $\sigma\colon G/G^p\to \CC^*$
such that $\sigma(x)=\omega$.
Now the term associated to $\tau$ in the first sum
equals the term associated to $\tau'=\tau\cdot\sigma$ in the second sum.

So we only have to consider terms of the form $x=\xi^p$ with $\xi\in G$.
The condition $\tau(x)=1$ is then independent of $\tau$ 
(since it is  equivalent to $\rho(\xi)=1$)
and the sum over all $\tau$ with $\psi\tau=\rho$ reduces
to a multiplication with the number of equivalence classes of such $\tau$.
By the Lemma we may assume that this number is nonzero.
Now $\tau_1$ and $\tau_2$ with $\psi\tau_1=\rho=\psi\tau_2$ are equivalent
iff $\tau_2=\tau_1^{1+mp^k}$ for some $m$ with $0\leq m<p$.
So this number is $1/p$ times the number of homomorphisms
$\sigma\colon G/G^p\to\CC^*$, hence equals $p^{r-1}$,
where $r$ denotes the rank of $G$.

On the other hand we have
\begin{equation*}
\psi^p b_\rho=\sum_{\rho\xi=1}\xi^p-\sum_{\rho\eta=\omega}\eta^p
\end{equation*}
Here the first sum is a certain factor times the sum over all $x\in G$
for which there exists $\xi\in G$ with $x=\xi^p$ 
and which satisfy $\tau(x)=1$ for any (and thus all) $\tau$ with $\psi\tau=\rho$.
The factor is the number of $\xi$ which satisfy these conditions,
which equals $p^r$.
\end{proof}

We write $h$ for the polynomial of degree $p-2$ given by
\begin{equation*}
h(t)=\frac{1}{1-t}\left(p-\frac{1-t^p}{1-t}\right)
\end{equation*}
\begin{proposition}
\begin{equation*}
\psi^p(b_1)=b_1+\sum_{\tau\not=1,\psi\tau=1}h(y_\tau)b_\tau
\end{equation*}
\end{proposition}

\begin{proof}
We have 
\begin{equation*}
b_\tau=\left(\sum_{\tau x=1} x\right)(1-y_\tau)
\end{equation*}
and thus
\begin{equation*}
h(y_\tau)b_\tau=\left(\sum_{\tau x=1}x\right)\left(p-\sum_{j=0}^{p-1}y_\tau^j\right)
=p\left(\sum_{\tau x=1}x\right)-\left(\sum_{x\in G}x\right)
\end{equation*}
We must take the sum of $\sum_{\tau x=1}x$
 over all equivalence classes of
$\tau\not=1$ with $\psi\tau=1$.
Interchange the sum over $x$ and the sum over $\tau$.
There are two cases:
\begin{itemize}
\item
If $x\in G^p$ then $\tau x=1$ for all $\tau$,
and we must simply count the number of equivalence class of $\tau$.
There are $p^r-1$ of them, with $p-1$ in each class.
\item
If $x\not\in G^p$ then the number of $\tau$ such that $\tau x=1$
is $p^{r-1}-1$, with again $p-1$ in each class.
\end{itemize}
So we get
\begin{equation*}
\begin{split}
\sum_{\tau\not=1,\psi\tau=1} h(y_\tau)b_\tau
=&p\left(\frac{p^r-1}{p-1}\sum_{x\in G^p}x+\frac{p^{r-1}-1}{p-1}\sum_{x\not\in G^p}x\right)\\
&-\left(\frac{p^r-1}{p-1}\sum_{x\in G^p}x+\frac{p^r-1}{p-1}\sum_{x\not\in G^p}x\right)\\
=&p^r\sum_{x\in G^p}x -\sum_{x\in G}x
=\sum_{\xi\in G}\xi^p-\sum_{x\in G}x=\psi^pb_1-b_1
\end{split}
\end{equation*}
\end{proof}

Now we consider the effect of Adams operations $\psi^q$ for 
primes $q\not=p$.
For any prime $q$ we write $f_q$ and $g_q$ for the polynomials given by
\begin{equation*}
f_q(t)=\frac{1-t^q}{1-t},\qquad
g_q(t)=\frac{(1-t)^{q-1}-f_q(t)}{q}
\end{equation*}

\begin{proposition}
If $\rho$ is of level $k>0$ then
\begin{equation*}
\psi^q(b_\rho)=f_q(y_\rho) b_\rho
\end{equation*}
and $\psi^q(b_1)=b_1$.
\end{proposition}

\begin{proof}
We have $ \psi^q(b_1)=\psi^q\left(\sum_{x\in G}x\right)
=\sum_{x\in G}x^q=\sum_{\xi\in G}\xi=b_1$ and
\begin{equation*}
\begin{split}
\psi^q(b_\rho)&=\psi^q\left(\left(\sum_{\rho x=1}x\right)\left(1-y_\rho\right)\right)
=\left(\sum_{\rho x=1}x^q\right)\left(1-y_\rho^q\right)\\
&=\left(\sum_{\rho\xi=1}\xi\right)(1-y_\rho)f_q(y_\rho)
=b_\rho f_q(y_\rho)
\end{split}
\end{equation*}
\end{proof}

\begin{corollary}
For the idempotents $e_\rho\in S_\rho$ one has
\begin{equation*}
\begin{split}
\psi^p(e_\rho)&=\sum_{\psi\tau=\rho} e_\rho\quad\text{ if }\rho\not=1,\\
\psi^p(e_1)&=e_1+\sum_{\tau\not=1,\psi\tau=1}e_\tau\\
\psi^q(e_\rho)&=e_\rho\quad\text{ if }q\not=p\\
\end{split}
\end{equation*}
\end{corollary}

Since $R$ has no $\ZZ$-torsion the Adams operations $\psi^\ell$ determine the
operations $\theta^\ell$ and we find
\begin{proposition}
\label{thetab}
If $\rho$ has level $k>0$  then
\begin{equation*}
\begin{split}
\theta^p(b_\rho)&=p^{(e-k)(p-1)-1}(1-y_\rho)^{p-1}b_\rho-\sum_{\psi\tau=\rho}b_\tau\quad\text{ if }k<e\\
\theta^p(b_\rho)&=g_p(y_\rho)b_\rho\quad\text{ if }k=e\\
\theta^q(b_\rho)&=\left(\frac{p^{(e-k)(q-1)}-1}{q}(1-y_\rho)^{q-1}+g_q(y_\rho)\right)b_\rho
\quad\text{ for }q\not=p\
\end{split}
\end{equation*}
Moreover
\begin{equation*}
\begin{split}
\theta^p(b_1)&=p^{e(p-1)-1}b_1-p^{-1}b_1-p^{-1}\sum_{\tau\not=1,\psi\tau=1}h(y_\tau)b_\tau\\
\theta^q(b_1)&=\frac{p^{e(q-1)}-1}{q}b_1
\quad\text{ for }q\not=p\\
\end{split}
\end{equation*}
\end{proposition}

\begin{proof}
This is just a matter of combining the last three Propositions
with the formula $\ell\theta^\ell(a)=a^\ell-\psi^\ell a$.
Note that $k=e$ can only happen if $G$ is cyclic,
in which case $y_\rho^p=1$,
which implies that 
$b_\rho^p=(1-y_\rho)^p=p(1-y_\rho)g_p(y_\rho)$.
\end{proof}

\begin{theorem}
\label{primconduc}
The $b_\rho$ with $\rho\not=1$ generate the $\lambda$-conductor ideal $I_\lambda$.
In other words $I_\lambda$ is the intersection of the augmentation ideal
and the ordinary conductor ideal $I$.
\end{theorem}

\begin{proof}
Write $J$ for the $R$-ideal generated by the $b_\rho$ with $\rho\not=1$.
From Proposition \ref{thetab} one reads of that $\theta^\ell(b_\rho)\in J$
for  $\rho\not=1$ and for every prim $\ell$.
From the identity
\begin{equation*}
\theta^\ell(ab)=\theta^\ell(a)b^\ell+\psi^\ell(a)\theta^\ell(b)
\end{equation*}
it then follows that $\theta^\ell(Rb_\rho)\subset J$ for $\rho\not=1$ and all $\ell$.
Finally from
\begin{equation*}
\theta^\ell(u+v)=\theta^\ell(u)+\theta^\ell(v)+\sum_{i=1}^{\ell-1}
\frac{1}{\ell}\binom{\ell}{i}u^iv^{\ell-i}
\end{equation*}
it follows that $\theta^\ell(J)\subset J$ for every $\ell$.
Since $J\subset I$ by Proposition \ref{conduc} this shows that $J\subset I_\lambda$.

Suppose that  $x\in I_\lambda$ and $x\not\in J$.
Then $x\in I$, so by Proposition \ref{conduc} there are $x_\rho\in R$ such that
$x=\sum_\rho x_\rho b_\rho$.
Since $\sum_{\rho\not=1}x_\rho b_\rho\in J\subset I_\lambda$ by the first half of the proof,
it follows that $x_1b_1\in I_\lambda$.
Since $gb_1=b_1$ for every $g\in G$ we may assume that $x_1\in\ZZ$.
Moreover $x_1\not=0$ which means that its $p$-valuation $v_p(x_1)$ is 
a natural number.
We may assume that $x$ is chosen in such a way that $v_p(x_1$) is minimal.
Now $I_\lambda$ must also contain
\begin{equation*}
\theta^p(x_1b_1)
=p^{-1}\left(x_1^p p^{e(p-1)}b_1-x_1(b_1+\sum_{\tau\not=1,\psi\tau=1}b_\tau)\right)
\end{equation*}
However the valuation of the coefficient of $b_1$ is
$v_p(x_1^p p^{e(p-1)}-x_1)-1=v_p(x_1)-1$,
in contradiction with the way $x$ was chosen.
Thus $I_\lambda\subset J$.
\end{proof}

\section {Direct products of relatively prime order}

Let $G_1$ be a group of order $n_1=p^e$,
and let $G_2$ a group of order $n_2=q^f$,
where $p$ and $q$ are different primes.
We write $R_1=\ZZ[G_1]$ and $R_2=\ZZ[G_2]$,
and denote their normal closures by $S_1$ and $S_2$ respectively.
Finally we write $I_1$ for the conductor of $S_1$ into $R_1$
and $I_2$ for the conductor of $S_2$ into $R_2$.
Since the $S_i$ are free abelian groups, the same is true
for the other additive groups involved, and we can view
$I_1\otimes I_2$ as a subgroup of $R_1\otimes I_2$ and of $R_1\otimes R_2$.
\begin{lemma}
\begin{equation*}
I_1\otimes I_2=(R_1\otimes I_2)\cap(I_1\otimes R_2)
\end{equation*}
\end{lemma}

\begin{proof}
There are $m_1,m_2\in\ZZ$ such that $m_1n_1+m_2n_2=1$.
If $x$ is an element of the left hand side then
$x\in R_1\otimes I_2$,
so $n_1x\in n_1R_1\otimes I_2\subset n_1S_1\otimes I_2\subset I_1\otimes I_2$
and therefore $m_1n_1x\in I_1\otimes I_2$.
Similarly $m_2n_2x\in I_1\otimes I_2$ and thus
$x=m_1n_1x+m_2n_2x\in I_1\otimes I_2$.
The other implication is obvious
\end{proof}

\begin{proposition}

The conductor $I$ of $S_1\otimes S_2$ into $R_1\otimes R_2$ is $I_1\otimes I_2$.
\end{proposition}

\begin{proof}
Suppose that $x\in I$, 
so that $x(S_1\otimes S_2)\subset R_1\otimes R_2$.
We write $x\in R_1\otimes R_2$ as $\sum x_g\otimes g$, 
where $g$ runs trough $G_2$.
For any $a\in S_1$ we have 
$\sum (x_ga)\otimes g=(\sum x_g\otimes g)(a\otimes 1)=x(a\otimes 1)\in R_1\otimes R_2$.
Therefore $x_ga\in I_1$ for any $a\in S_1$,
which means that $a_g\in I_1$ for all $g\in G_1$.
Thus $x\in I_1\otimes R_2$.
Similarly $x\in R_1\otimes I_2$.
Thus $x\in I_1\otimes I_2$ by the Lemma.
The other inclusion is obvious.
\end{proof}

We show now that for the $\lambda$-conductor
 a similar theorem holds:
 
 \begin{theorem}
 The $\lambda$-conductor $I_\lambda$ of $S_1\otimes S_2$ into $R_1\otimes R_2$
 is the tensor product of the $\lambda$-conductors $I_{1\lambda}$ of $S_1$ into $R_1$
 and $I_{2\lambda}$ of $S_2$ into $R_2$.
 \end{theorem}
 
\begin{proof}
The $\lambda$-conductor $I_\lambda$  is a subset of the classical conductor $I$ ,
which is $I_1\otimes I_2$.
However $I_1$ is the direct sum $\ZZ b_1\oplus I_{1\lambda}$ by theorem \ref{primconduc}.
and similarly for $I_2$.
Thus any $x\in I_\lambda$ can uniquely be written as
\begin{equation*}
x=x_0(b_1\otimes b_1)\oplus(x_1\otimes b_1)\oplus(b_1\otimes x_2)\oplus y
\end{equation*}
with $x_0\in\ZZ$, 
$x_1\in I_{1\lambda}$, 
$x_2\in I_{2\lambda}$, 
$y\in I_{1\lambda}\otimes I_{2\lambda}$.
Since $I_\lambda$ is an ideal of $S_1\otimes S_2$, each of these four
summands must be in $I_\lambda$.

Therefore we consider the intersection of $I_\lambda$ with $b_1\otimes I_{2\lambda}$.
Suppose that $a$ is an element of this intersection, 
say $a=b_1\otimes x$ with $x\in I_{2\lambda}$.
Then $\theta^p(a)\in I_\lambda$ too.
We have
\begin{equation*}
\theta^p(a)
=p^{-1}(a^p-\psi^pa)
=p^{-1}(p^{e(p-1)}b_1\otimes x^p-(b_1+\sum_{\tau\not=1,\psi\tau=1}b_\tau)\otimes\psi^px)
\end{equation*}
and thus $p^{e(p-1)-1}b_1\otimes x^p-p^{-1}b_1\otimes\psi^px $ 
should be in $I_\lambda$.
Now the first term is a multiple of $a$ and thus in $I_\lambda$.
So the other term $p^{-1}b_1\otimes\psi^px$ is in the aforementioned intersection.
Since $\psi^p$ is an automorphism (of finite order) of $R_2$
this shows that the intersection is $p$-divisible.
Since the intersection is a finitely generated abelian group 
this can only happen if it vanishes.

The same argument applies to the first and second summand of $x$.
Thus $x=y\in I_{1\lambda}\otimes I_{2\lambda}$ and we have shown that
$I_\lambda\subset I_{1\lambda}\otimes I_{2\lambda}$.
The other inclusion is obvious.
\end{proof}



\end{document}